# Benders decomposition approach to solve the capacitated facility location problem

Ali Akbar Sadat Asl[*,1], Ali Rouhani[2]

Department of Industrial Engineering, Amirkabir University of Technology, Tehran, Iran[1]
Department of Computer Science and Operations Research, University of Montreal, Quebec, Canada[2]

**Abstract**

Facility Location (FL) problems as one of the most important problems in operations research aim to determine the location of a set of facilities in a way that the total costs, including costs of opening facilities and transportation costs, are minimized. This study addresses an FL problem in which the capacity of each facility is limited. Because this problem is in the category of np-hard problems, we use the Benders Decomposition (BD) approach to efficiently solve the FL problem. In this paper, we implement the classic BD algorithm and some accelerating BD methods such as Pareto-optimality cut and L-shaped decomposition methods. Furthermore, we propose and implement the hybrid Pareto-L-shaped (PL) method, and evaluate the performance of the implemented algorithms. The results show that the L-shaped decomposition outperforms the other algorithms in terms of time and the number of iteration, while the classic BD converges slowly especially on large scales.

**Keywords:** Facility Location, Benders Decomposition, Pareto-optimality Cut, L-shaped Decomposition

## 1. Introduction

FL problem is one of the most important problems in operations research and resource allocation that has attracted much attention in recent decades. This problem has many applications in the real world such as the construction of hospitals and health centers, construction of fire stations, establishing warehouses and factories, and finding suitable places for amusement parks [1]. The FL problem is a combinatorial optimization problem to determine the optimal locations of a set of facilities (warehouses, plants, machines, hospitals, fire stations, etc.) and assign customers to these facilities in such a way that the total cost is minimized [2].

In this problem, two types of costs are considered. The first cost is related to the setup cost (facility cost) that occurs when a facility is opened. The second one is a transportation cost which occurs when a customer is allocated to an opened facility. For example, consider we want to construct some warehouses. The total costs are comprised of the cost of opening the warehouses and the transportation cost associated with shipping products from warehouses to customers. FL problems are classified into two types: Uncapacitated Facility Location Problem (UFLP) and Capacitated Facility Location Problem (CFLP). If an arbitrary number of customers can be connected to a facility, the problem is called UFLP. In contrast, if each facility has a limit on the number of customers it can serve, it becomes a CFLP [3].

Both UFLP and CFLP are np-hard problems. The FL problems have been extensively studied, and to solve them, a lot of exact and heuristic algorithms have been developed in the past decades [4]. Some exact algorithms such as branch-and-bound algorithm [5], Lagrangian-based algorithm [6], cutting plane method [7], Dantzig–Wolfe decomposition [8], and BD algorithm [9] have been developed to solve facility location problems. Another line of research lies in developing heuristic and approximation algorithms, such as genetic algorithm [10], tabu search algorithm [11], and particle swarm optimization [12].

In this study, we use the BD approach to solve the CFLP. A great advantage of BD is that it converges straight to the optimal of the Mixed Integer Linear Program (MILP) rather than to a relaxation of the problem, as Dantzig–Wolfe decomposition and Lagrangian relaxation do. Therefore, BD does not need to be embedded within a branch and bound framework [13]. BD algorithm is one of the most widely used exact algorithms in dealing with difficult optimization problems. BD partitions the problem into multiple sub-problems instead of considering all decision variables and constraints of a large-scale problem simultaneously by which we can reduce the complexity of the problem, and the optimal solution is guaranteed yet [14].

Although BD is one of the most important exact algorithms in dealing with difficult optimization problems, slow convergence is the main drawback of this algorithm. To deal with this issue, many researchers have worked on accelerating methods. Pareto-optimality cut and L-shaped methods are the most common approaches which can be used in many problems. In this paper, in addition to the implementation of classic BD, Pareto optimality cut, and L-shaped algorithms, we combine Pareto and L-shaped methods and propose the hybrid PL method. To evaluate the performance of these algorithms, several cases are solved using these algorithms, and the results, are compared.

The structure of this paper is as follows: the CFLP is addressed in section 2. The procedure of the BD algorithm is described in section 3. In section 4, the accelerating BD methods are presented. The performance of the implemented algorithms is evaluated in section 5. Finally, discussion and conclusions are considered in section 6.

## 2. Capacitated Facility Location Problem (CFLP)

CFLP addresses the problem of locating a new set of facilities in a way that the total costs comprised of transportation cost and cost of opening the facilities are minimized. CFLP is a variant of the FL problem in which we consider capacities for the facilities. In other words, with the inclusion of the capacities, an open facility may not be able to serve an arbitrary number of customer demands. Table 1 shows sets, parameters, and decision variables in the CFLP.

**Table 1. Sets, parameters, and decision variables of CFLP**

| Symbol | Description |
|---|---|
| Sets | |
| $M$ | Set of demand nodes (customers) $i \in M$ |
| $N$ | Set of facilities $j \in N$ |
| Parameters | |
| $f_j$ | Fixed cost of opening a facility on node $j$ |
| $s_j$ | Capacity of facility $j$ |
| $d_i$ | Demand of customer $i$ |
| $c_{ij}$ | Transportation cost between customer $i$ and facility $j$ |
| Decision Variables | |
| $y_j$ | If a facility at node $j$ is opened 1, otherwise 0 |
| $x_{ij}$ | Fraction of the demand of customer $i$ served by facility $j$ |

The mathematical model of CFLP based on the abovementioned symbols is as follows:

$$\min_{x,y} \sum_{j=1}^{n} f_j y_j + \sum_{i=1}^{m}\sum_{j=1}^{n} d_i c_{ij} x_{ij} \quad (1)$$

$$\sum_{j=1}^{n} x_{ij} = 1; \quad \forall i \in M \quad (2)$$

$$x_{ij} \leq y_j; \quad \forall i \in M, \forall j \in N \quad (3)$$

$$\sum_{i=1}^{m} d_i x_{ij} \leq s_j y_j; \quad \forall j \in N \quad (4)$$

$$x_{ij} \geq 0; \quad \forall i \in M, \forall j \in N \quad (5)$$

$$y_j \in \{0, 1\}; \quad \forall j \in N \quad (6)$$

In the above model, the objective function (1) minimizes the total costs comprised of transportation cost and cost of opening the facilities. Constraints (2) imply that the demand for each customer should be fully satisfied. Constraints (3) state, before allocating some nodes to the node j, we should have a facility at node j. Constraints (4) say that the capacity of facilities is limited; in other words, the total demands of all customers that should be satisfied at node j, should not exceed the capacity of the facility at node j, in case that node j is a facility node. Finally, Constraints (5) and (6) state the condition of decision variables of the problem.

To illustrate this problem, consider that we have some nodes. These nodes are called demand nodes or customers. In this problem, we should decide about two variables. The first one is related to the location of facilities. Suppose that we want to construct four warehouses. In this case, our facilities are warehouses. The first thing that should be determined is that in which of these nodes we should construct the warehouses. For example, assume that we construct the facilities in nodes 3, 4, 7, and 9. Now, the demands of nodes 3, 4, 7, 9 are satisfied because they have a facility. The second thing that should be determined is that how we should allocate these nodes in which there isn't a facility to the facility nodes. For example, here, sixty percent of the demand of node 10 is satisfied by the facility in node 7, and the rest of that is satisfied by the facility established in node 9, and node 1 is fully satisfied by the facility opened in node 3.

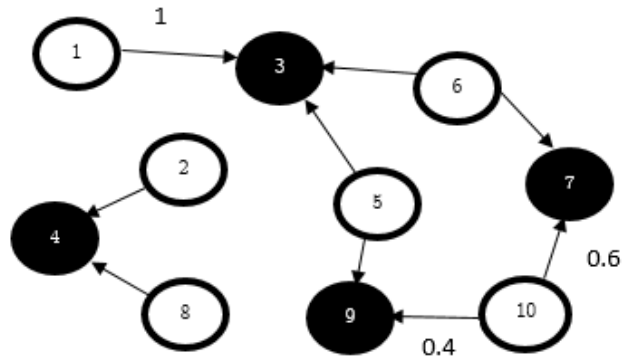

**Figure 1. An example of CFLP**

## 3. Benders Decomposition

BD algorithm is a partitioning method to solve large-scale mixed-integer programming (MIP) problems with complicated variables. Before going to the details of the algorithm, we define the term "complicated variable". Consider we have an optimization problem defined as follows:

$$\min_{x,y} c^T x + d^T y$$
$$A_1 x_1 + E_1 y \geq b_1 \quad (7)$$
$$A_2 x_2 + E_2 y \geq b_2$$

In model (7), variable $y$ is a complicated variable because it prevents us from partitioning the model into some sub-problems. If fixing this variable to a constant value, we could partition this problem into these two sub-problems ($\bar{y}$ is a fixed value for $y$):

$$\min_{x_1} c^T x_1 \qquad\qquad \min_{x_2} c^T x_2$$
$$A_1 x_1 \geq b_1 - E_1 \bar{y} \qquad A_2 x_2 \geq b_2 - E_2 \bar{y}$$

BD is based on a relaxation algorithm that splits the original problem into two smaller programs: a master problem and a set of sub-problems, taking advantage of dual decomposition structure. The master problem is a relaxed version of the original problem with the sets of integer variables, known as complicated variables, and its associated constraints. The sub-problems are dual versions of what would be the original problem if the values of the integer variables were temporarily fixed by the master problem. In other words, we solve the master problem and obtain a value for the complicated variables. Then, we give the results of solving the master problem to the sub-problem, and the sub-problem would be an LP model instead of a MIP problem as the complicated variables have been fixed to the obtained values from the master problem. BD algorithm solves the master problem iteratively followed by the implied sub-problems. At each iteration, a new valid inequality (called a Benders cut) derived from the sub-problems is added to the master program. The algorithm converges to an optimal solution for the original mixed-integer program if one exists [15]. Figure 2 shows the way the BD algorithm works.

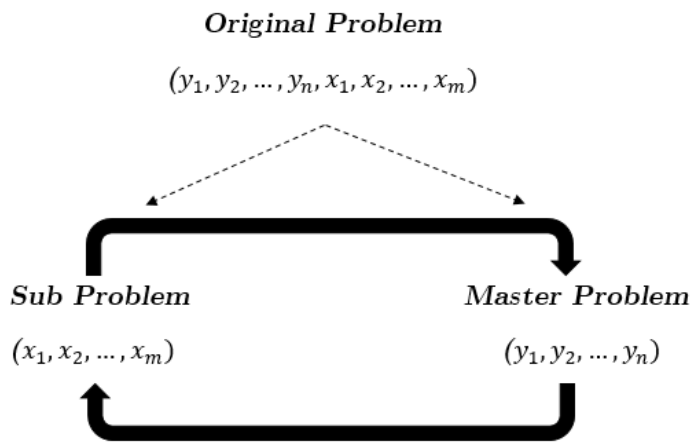

**Figure 2. An overview of the BD algorithm**

Now, we explain above procedure in an MIP model defined as follows:

$$\min_{x,y} c^T x + f^T y$$
$$Ax + By \geq b$$
$$y \in Y$$
$$x \geq 0$$
(8)

If $y$ is fixed to a feasible integer configuration, the resulting model to solve is [16]:

$$\min_{x} c^T x$$
$$Ax \geq b - B\bar{y}$$
$$x \geq 0$$
(9)

Therefore, the complete minimization problem can be written as:

$$\min_{y \in Y}[f^T y + \min\{c^T x | Ax \geq b - By\}]$$
(10)

The dual of the inner LP problem is:

$$\max_{u}(b - B\bar{y})^T u$$
$$A^T u \leq c$$
$$u \geq 0$$
(11)

The BD algorithm can be stated as [16]:

{*initialization*}

$y :=$ *initial feasible integer solution*

$LB := -\infty$  (*LB is a lowerbound for the objective function of the original problem*)

$UB := +\infty$  (*UB is an upperbound for the objective function of the original problem*)

**while** $UB - LB > \epsilon$ **do**

   {*solve subproblem*}

   $\min\{f^T \bar{y} + (b - B\bar{y})^T u | A^T u \leq c, u \geq 0\}$

   **if** *Unbounded* **then**

      *Get unbounded ray* $\bar{u}$

      *Add cut* $(b - B\bar{y})^T \bar{u} \leq 0$ *to master problem*

   **else**

      *Get extreme point* $\bar{u}$

      *Add cut* $z \geq f^T y + (b - By)^T \bar{u}$ *to master problem*

      $UB := \min\{UB, f^T \bar{y} + (b - B\bar{y})^T \bar{u}\}$

***end if***

{*solve master problem*}

$min\{z|cuts, y \in Y\}$

$LB := z$

***end while***

The sub-problem is a dual LP problem, and the master problem is a pure IP problem (no continuous variables are involved). BD for MIP is of special interest when the sub-problem and the relaxed master problem are easy to solve, while the original problem is not [16].

Now, consider the CFLP discussed earlier. We can rewrite the mathematical model as:

$$\min_{x,y} \sum_{j=1}^{n} f_j y_j + \sum_{i=1}^{m} \sum_{j=1}^{n} d_i c_{ij} x_{ij} \tag{12}$$

$$\sum_{j=1}^{n} x_{ij} \geq 1; \quad \forall i \in M \tag{13}$$

$$-x_{ij} \geq -y_j; \quad \forall i \in M, \forall j \in N \tag{14}$$

$$-\sum_{i=1}^{m} d_i x_{ij} \geq -s_j y_j; \quad \forall j \in N \tag{15}$$

$$x_{ij} \geq 0; \quad \forall i \in M, \forall j \in N \tag{16}$$

$$y_j \in \{0,1\}; \quad \forall j \in N \tag{17}$$

The BD sub-problem can be stated as:

$$\max_{u,v,w} \sum_i u_i - \sum_i \sum_j y_j v_{ij} - \sum_j s_j y_j w_j$$

$$u_i - v_{ij} - d_i w_j \leq d_i c_{ij} \tag{18}$$

$$u_i, v_{ij}, w_j \geq 0$$

where $u_i, v_{ij}, w_j$ are dual variables related to constraints (13), (14), and (15), respectively.

The BD relaxed master problem can be written as:

$$\min_y z$$

$$z \geq \sum_j f_j y_j + \sum_i \overline{u_i} - \sum_i \sum_j y_j \overline{v_{ij}} - \sum_j s_j y_j \overline{w_j} \tag{19}$$

$$\sum_i \overline{u_i} - \sum_i \sum_j y_j \overline{v_{ij}} - \sum_j s_j y_j \overline{w_j} \leq 0$$

$$y_j \in \{0,1\}$$

# 4. Accelerating Benders Decomposition

Although BD is one of the most important exact algorithms in dealing with difficult optimization problems, slow convergence is the main drawback of this algorithm. To deal with this issue, many researchers have worked on accelerating methods. Pareto-optimality cut and L-shaped methods are the most common approaches which can be used in many problems. These algorithms follow a procedure similar to the classic BD, but in the Pareto-optimality cut, we select the strongest cuts in each iteration. In the L-shaped method, we partition the sub-problem into, for example, L sub-problems. So, instead of having only one cut in each iteration, we can generate L cuts.

## 4.1. Pareto-Optimality Cut

A cut is called Pareto-optimal if no other cut dominates it. A Pareto-optimal cut only exists when the dual of the Benders sub-problem has several optimal solutions, and it is the strongest cut among all the alternative Benders cuts in the same iteration. Consequently, Pareto-optimal cuts can effectively improve the performance of the BD algorithm [17].

In other words, given the two sets of dual variable solutions to a sub-problem, $u^1$ and $u^2$, the cut $z \geq f(u^1) + g(u^1)y$ is said to dominate or is considered stronger than the cut $z \geq f(u^2) + g(u^2)y$ if $f(u^1) + g(u^1)y \geq f(u^2) + g(u^2)y$ for all $y$ in the polytope $P$, with strict inequality for at least one $y \in P$. A dual variable solution $u^1$ is said to dominate $u^2$ if the associated cut is stronger, and $u^1$ is called Pareto-optimal if its corresponding cut is Pareto-optimal [18].

As mentioned earlier, the Pareto-optimal cut only exists in which the dual of the Benders sub-problem has several optimal solutions. The following theorem provides a method for choosing from among the alternate optimal solutions to the sub-problem to generate Pareto-optimal cuts. Let $Y^c$ the relative interior of the convex hull of the $y$-space of the master problem. Consider that $U(\tilde{y})$ is the set of optimal solutions to the following optimization problem:

$$\max\{f(u) + g(u)\tilde{y}\}$$
$$u \in U \tag{20}$$

Then, $u^o$, the optimal solution of the following model, will result in generating the Pareto optimization cut.

$$\max\{f(u) + g(u)y^o\}$$
$$u \in U(\tilde{y}) \tag{21}$$

where $y^o$ is a core point contained in $Y^c$ (for more details please refer to [18]).

## 4.2. L-shaped Decomposition Method

The L-shaped method is a derivation of the BD algorithm in which the sub-problem is decomposed to some sub-problems. Therefore, instead of having only one cut in each iteration, we have a family of cuts contributing to faster convergence. One of the successful applications of this method is in scenario-based stochastic programming problems. In these problems, we have k scenarios and based on the scenarios, the sub-problem is partitioned into k sub-problems [19].

In CFLP, we have $m$ customers. If decomposing the sub-problem based on each customer, we will have $m$ sub-problems, each of which has a simpler solution process (see figure 3).

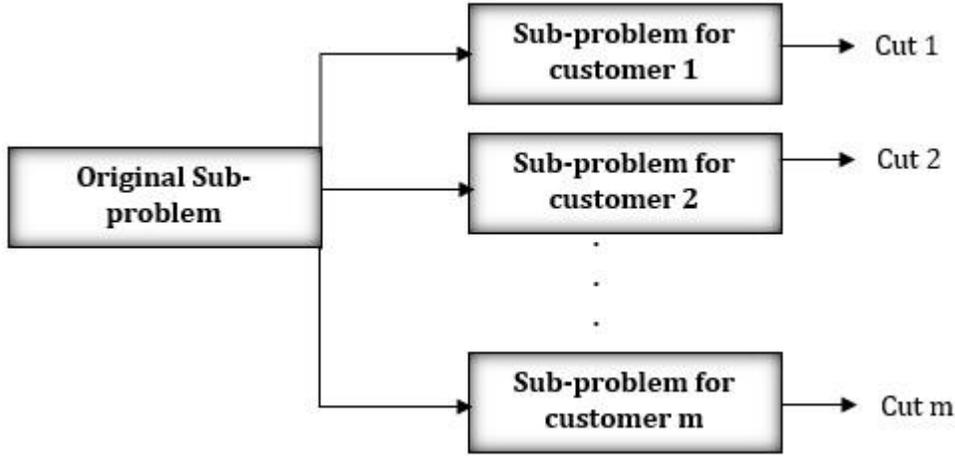

**Figure 3. L-shaped decomposition**

### 4.3. Hybrid Pareto-optimality Cut and L-shaped Decomposition Method

In this paper, we combined the Pareto-optimality cut and L-shaped decomposition to generate the hybrid PL algorithm. Therefore, in each iteration of the algorithm, in addition to the selection of the strongest cut among all the alternative Benders cuts, we also partition the original sub-problem into $m$ sub-problems. As we obtain the strongest cuts and decompose the sub-problem, more computations are needed in this hybrid approach. However, we expect that PL can converge to the optimal solution in a fewer number of iterations as we make use of the advantages of Pareto-optimality cut and L-shaped decomposition at the same time.

## 5. Computational Results

In this study, we implemented the classic BD, Pareto-optimality cut, L-shaped decomposition method, and the hybrid PL algorithm. To evaluate the performance of these four algorithms, five cases have been generated and solved using these algorithms implemented on the General Algebraic Modeling System (GAMS). Table 2 shows these cases.

**Table 2. Five cases solved by the implemented algorithms**

|  | Case 1 | Case 2 | Case 3 | Case 4 | Case 5 |
|---|---|---|---|---|---|
| Number of Customers | 5 | 10 | 50 | 70 | 70 |
| Number of Facilities | 2 | 4 | 20 | 20 | 30 |

In each case, the values of the parameter $c_{ij}$ (transportation cost between customer i and facility j) are achieved by the uniform distribution in (50, 100). The values of $f_j$ (fixed cost of opening a facility on node j) are obtained by the uniform distribution in (1000, 1500). Furthermore, the values of $d_i$ (demand of customer $i$) and $s_j$ (capacity of facility $j$) are uniformly obtained from intervals (50, 100) and (2000, 2500), respectively. All of the algorithms ran on a personal computer with a 2.53 GHz processor and 6 GB of RAM.

Table 3 compares the performance of the implemented algorithms in terms of time and the number of iteration. It is noteworthy to mention that all these algorithms converged to the optimal solution. The values of optimal costs can be seen in Table 3. In terms of the number of iteration, the hybrid PL and L-shaped decomposition methods outperform the other algorithms. As can be seen, the classic BD requires a large number of iterations to converge to the optimal solution.

In terms of time, although the Pareto-optimality cut outperforms the other algorithms in cases 1, 2, and 3, it turned out that the L-shaped decomposition can perform better in larger scales (cases 4 and 5). In case 5, while the L-shaped method needs only about 13 seconds to converge to the optimal solution, the classic BD entails about 12234 seconds, which is a long time compared to the other algorithms.

Overall, considering both performance measures (time and the number of iteration), we notice that the L-shaped decomposition method outperforms the other algorithms for the CFLP.

Table 3. Comparison of implemented algorithms

| Criteria | Algorithms | Case 1 | Case 2 | Case 3 | Case 4 | Case 5 |
|---|---|---|---|---|---|---|
| Iteration | BD | 5 | 5 | 117 | 231 | 518 |
| | Pareto | 1 | 1 | 30 | 46 | 107 |
| | L-shaped | 1 | 3 | 7 | 5 | 6 |
| | Hybrid PL | 1 | 1 | 6 | 5 | 7 |
| Time | BD | 0.76 | 0.47 | 33.34 | 574.83 | 12234.88 |
| | Pareto | 0.05 | 0.06 | 8.59 | 17.23 | 74.312 |
| | L-shaped | 0.14 | 0.81 | 11.53 | 10.72 | 13.34 |
| | Hybrid PL | 0.14 | 0.44 | 21.40 | 20.08 | 16.49 |
| Optimal Cost | | 26304.99 | 47276.46 | 210638.24 | 287097.83 | 279559.45 |

## 6. Discussion and Conclusion

In this study, we address the CFLP, which is an NP-hard problem. To efficiently solve this problem we use the BD approach by which the original problem is partitioned into multiple smaller problems. As the difficulty of optimization problems is commensurate with the number of variables and constraints, iteratively solving these smaller problems could be more efficient than solving a single large problem. However, the BD algorithm converges to the optimal solution slowly, and to deal with this issue, some accelerating methods can be used. Pareto-optimality cut and L-shaped methods are the most common approaches which can be used in many problems.

In this paper, we implement the classic BD, Pareto-optimality cut, L-shaped decomposition, and the hybrid PL methods for the CFLP. Some cases have been generated randomly and solved by these algorithms. The results show the classic BD is outperformed by the accelerating methods in terms of time and the number of iteration. Comparing the implemented accelerating methods, the hybrid PL and L-shaped decomposition methods perform better than the other algorithms in terms of the number of iteration. In terms of time, although the Pareto-optimality cut performs quickly in small scales, this algorithm is outperformed by the L-shaped method in larger scales.

**References**

[1] Sadatasl, Ali Akbar, Mohammad Hossein Fazel Zarandi, and Abolfazl Sadeghi. "A combined facility location and network design model with multi-type of capacitated links and backup facility and non-deterministic demand by fuzzy logic." 2016 Annual Conference of the North American Fuzzy Information Processing Society (NAFIPS). IEEE, 2016.

[2] Sadat Asl, A. A., et al. "A fuzzy capacitated facility location-network design model: A hybrid firefly and invasive weed optimization (FIWO) solution." Iranian Journal of Fuzzy Systems 17.2 (2020): 79-95.

[3] Wu, Ling-Yun, Xiang-Sun Zhang, and Ju-Liang Zhang. "Capacitated facility location problem with general setup cost." Computers & Operations Research 33.5 (2006): 1226-1241.

[4] Li, Yan, et al. "Optimal physician assignment and patient demand allocation in an outpatient care network." Computers & Operations Research 72 (2016): 107-117.

[5] Tcha, Dong-wan, and Bum-il Lee. "A branch-and-bound algorithm for the multi-level uncapacitated facility location problem." European Journal of Operational Research 18.1 (1984): 35-43.

[6] Beasley, John E. "An algorithm for solving large capacitated warehouse location problems." European Journal of Operational Research 33.3 (1988): 314-325.

[7] Avella, Pasquale, and Maurizio Boccia. "A cutting plane algorithm for the capacitated facility location problem." Computational Optimization and Applications 43.1 (2009): 39-65.

[8] Wu, Tao, et al. "Dantzig-Wolfe decomposition for the facility location and production planning problem." Computers & Operations Research 124 (2020): 105068.

[9] Fischetti, Matteo, Ivana Ljubić, and Markus Sinnl. "Benders decomposition without separability: A computational study for capacitated facility location problems." European Journal of Operational Research 253.3 (2016): 557-569.

[10] Fernandes, Diogo RM, et al. "A simple and effective genetic algorithm for the two-stage capacitated facility location problem." Computers & Industrial Engineering 75 (2014): 200-208.

[11] Sun, Minghe. "A tabu search heuristic procedure for the capacitated facility location problem." Journal of Heuristics 18.1 (2012): 91-118.

[12] Chu, Xianghua, et al. "An orthogonal-design hybrid particle swarm optimiser with application to capacitated facility location problem." International Journal of Bio-Inspired Computation 8.5 (2016): 268-285.

[13] Ibrahim, Ahmed, et al. "Bender's Decomposition for Optimization Design Problems in Communication Networks." IEEE Network 34.3 (2019): 232-239.

[14] Taskın, Z. Caner. "Benders decomposition." Wiley Encyclopedia of Operations Research and Management Science. John Wiley & Sons, Malden (MA) (2010).

[15] de Camargo, Ricardo Saraiva, Gilberto de Miranda Jr, and Henrique Pacca L. Luna. "Benders decomposition for hub location problems with economies of scale." Transportation Science 43.1 (2009): 86-97.

[16] Kalvelagen, Erwin. "Benders decomposition with GAMS." Amsterdam Optimization Modeling Group: Washington, DC, USA (2002).

[17] Tang, Lixin, Wei Jiang, and Georgios KD Saharidis. "An improved Benders decomposition algorithm for the logistics facility location problem with capacity expansions." Annals of operations research 210.1 (2013): 165-190.

[18] Watson, Frederick R., and David F. Rogers. "Pareto-Optimality of the Balinski Cut for the Uncapacitated Facility Location Problem." (2006).

[19] Shiina, Takayuki. "L-shaped decomposition method for multi-stage stochastic concentrator location problem." Journal of the Operations Research Society of Japan 43.2 (2000): 317-332.